\documentclass[a4paper,11pt]{amsart}
\usepackage{amsmath,amssymb,amsfonts,latexsym}

\newtheorem{theorem}{Theorem}[section]

\input{tcilatex}

\begin{document}
\title[On a discontinuous retarded differential operator]{Spectrum, trace
and oscillation of a Sturm-Liouville type retarded differential operator
with interface conditions}
\author{Erdo\u{g}an \c{S}en}
\address{Department of Mathematics, Nam\i k Kemal University, 59030, Tekirda%
\u{g}, Turkey}
\email{erdogan.math@gmail.com}

\begin{abstract}
In this study, a formula for regularized sums of eigenvalues of a
Sturm-Liouville problem with retarded argument at the point of discontinuity
is obtained. Moreover, oscillation properties of the related problem is
investigated.


\vspace{2mm}

\end{abstract}

\maketitle

%%%%%%%%%%%%%%%%%%%%%%%%%%%%%%%%%%%%%%%%%%%%%%%%%%%%%%%%%%%%%%%%%%%

%%%%%%%%%%%%%%%%%%%%%%%%%%%%%%%%%%%%%%%%%%%%%%%%%%%%%%%%%%%%%%%%%%%

\section{Introduction}

%%%%%%%%%%%%%%%%%%%%%%%%%%%%%%%%%%%%%%%%%%%%%%%%%%%%%%%%%%%%%%%%%%%

In this paper, we consider the boundary value problem for the differential
equation with retarded argument%
\begin{equation}
p(x)y^{\prime \prime }(x)+q(x)y(x-\Delta (x))+\lambda ^{2}y(x)=0,  \tag{1}
\end{equation}%
on $\left[ 0,\frac{\pi }{2}\right) \cup \left( \frac{\pi }{2},\pi \right] ,$
with boundary conditions%
\begin{equation}
a_{1}y(0)+a_{2}y^{\prime }\left( 0\right) =0,  \tag{2}
\end{equation}%
\begin{equation}
y^{\prime }(\pi )+dy\left( \pi \right) =0  \tag{3}
\end{equation}%
and interface conditions%
\begin{equation}
\gamma _{1}y\left( \frac{\pi }{2}-0\right) -\delta _{1}y(\frac{\pi }{2}+0)=0,
\tag{4}
\end{equation}%
\begin{equation}
\gamma _{2}y^{\prime }\left( \frac{\pi }{2}-0\right) -\delta _{2}y^{\prime }(%
\frac{\pi }{2}+0)=0,  \tag{5}
\end{equation}%
where $p(x)=p_{1}^{2}$ for $x\in \left[ 0,\frac{\pi }{2}\right) $ and $%
p(x)=p_{2}^{2}$ for $x\in \left( \frac{\pi }{2},\pi \right] $; the
real-valued function $q\left( x\right) $ is continuous in $\left[ 0,\frac{%
\pi }{2}\right) \cup \left( \frac{\pi }{2},\pi \right] $ and has finite
limits $q\left( \frac{\pi }{2}\pm 0\right) =\lim_{x\rightarrow \frac{\pi }{2}%
\pm 0}q\left( x\right) ,$ the real-valued function $\Delta (x)\geq 0$ is
continuous in $\left[ 0,\frac{\pi }{2}\right) \cup \left( \frac{\pi }{2},\pi %
\right] $ has finite limits $\Delta \left( \frac{\pi }{2}\pm 0\right)
=\lim_{x\rightarrow \frac{\pi }{2}\pm 0}\Delta \left( x\right) $, \textit{if}%
$\ x\in \left[ 0,\frac{\pi }{2}\right) $ then $x-\Delta (x)\geq 0$; \textit{%
if} $x\in \left( \frac{\pi }{2},\pi \right] $ then $x-\Delta (x)\geq \frac{%
\pi }{2}$;$\lambda $ is a spectral parameter; $p_{i},a_{i},d,\gamma
_{i},\delta _{i}$ $\left( i=1,2\right) $ are arbitrary real numbers such
that $p_{i}a_{i}d\neq 0$ $(i=1,2).$

Differential equations with retarded argument appeared as far back as the
eighteenth century in connection with the solution of the problem of Euler
on the investigation of the general form of curves similar to their own
evolutes. Differential equations with retarded argument arise in many
application of mathematical modelling: for example, combustion in the
chamber of a liquid propellant rocket engine [1,2] and vibrations of the
hammer in an electromagnetic circuit breaker [3,4].

Boundary value problems with interface conditions arise in varied assortment
of physical transfer problems (see [5]). Also, some problems with interface
conditions arise in thermal conduction problems for a thin laminated plate
(i.e., a plate composed by materials with different characteristics piled in
the thickness). In this class of problems, interface conditions across the
interfaces should be added since the plate is laminated. The study of the
structure of the solution in the matching region of the layer with the basis
solution in the plate leads to consideration of an eigenvalue problem for a
second-order differential operator with piecewise continuous coefficients
and interface conditions [6].

The asymptotic formulas for eigenvalues of boundary value problems with
retarded argument obtained in [7-17]. In [17], the principal term of
asymptotic distribution of eigenvalues of the problem (1)-(5) was obtained
up to $O\left( \frac{1}{n}\right) $. But we need sharper asymptotic
formulas. Therefore, we improve this formula up to $O\left( \frac{1}{n^{3}}%
\right) $.

The theory of regularized traces of Sturm-Liouville operators stems from the
paper [18] of Gelfand and Levitan. Trace formulas for the Sturm-Liouville
equation with a complex valued potential and with two point boundary
conditions obtained in [19]. Regularized trace of the Sturm-Liouville
problem with irregular boundary conditions investigated in [20]. A
regularized trace formula for the matrix Sturm-Liouville operator found in
[21]. The first regularized traces of boundary value problems with unbounded
operator coefficient in a finite interval obtained in [22-24]. The second
regularized trace of a differential operator with bounded operator
coefficient and with the mixed boundary conditions can be found in [25]. For
a comprehensive review of traces of Schr\"{o}dinger operators, the
interested reader is referred to [26].

As mentioned above partly, the literature is about the regularized traces of
classic type differential operators is so rich and diverse. For a more
comprehensive list, one can refer to the survey paper [27] and the paper
[28]. However, there are only a few works on regularized traces and
oscillation properties for differential operators with retarded argument. M.
Pikula in [8] obtained trace formula of first order:

if $\tau \geq \pi $%
\begin{equation*}
\sum_{n=1}^{\infty }\left( \lambda _{n}\left( \tau \right) -n^{2}\right)
=q\left( 0\right) \varphi \left( -\tau \right) -\frac{h^{2}+H^{2}}{2}
\end{equation*}%
and if $\tau \leq \pi $%
\begin{eqnarray*}
&&\sum_{n=1}^{\infty }\left[ \lambda _{n}\left( \tau \right) -n^{2}-\frac{2}{%
\pi }\left( h+H+\frac{\cos n\tau }{2}\dint\limits_{\tau }^{\pi }q\left(
t\right) dt\right) \right] \\
&=&q\left( 0\right) \varphi \left( -\tau \right) -\frac{h^{2}+H^{2}}{2}+%
\frac{h+H}{2\pi }\dint\limits_{\tau }^{\pi }q\left( t\right) dt+\left[ \frac{%
q\left( \tau \right) +q\left( \pi \right) }{4}-q\left( \tau \right) \right] 
\frac{\pi -\tau }{\pi } \\
&&+\left[ \frac{b}{4}-\left( \frac{1}{\sqrt{8}}\dint\limits_{\tau }^{\pi
}q\left( t\right) dt\right) ^{2}\right] \frac{\pi -2\tau }{\pi }
\end{eqnarray*}%
for the boundary-value problem of second order with retarded argument:%
\begin{gather*}
-y^{\prime \prime }+q(x)y(x-\tau )=\lambda y, \\
y^{\prime }(0)-hy(0)=y^{\prime }(\pi )+Hy(\pi )=0, \\
y(x-\tau )=y(0)\varphi (x-\tau ),\quad x\leq \tau ,\;\varphi (0)=1.
\end{gather*}

C-F. Yang in [11] obtained formula of the first regularized trace,
oscillations of the eigenfunctions and the solutions of inverse nodal
problem for discontinuous boundary value problems with retarded argument and
with interface conditions at the one point of discontinuity. F. Hira in [14]
obtained a formula for regularized sums of eigenvalues for a Sturm-Liouville
problem with retarded argument at the one point of discontinuity which
contains a spectral parameter in the boundary conditions and as a most
recent study in this topic, \c{S}en studied regularized trace formula and
oscillation of eigenfunctions of a Sturm-Liouville operator with retarded
argument at two points of discontinuity [16].

The goals of this article are to calculate the regularized trace and to find
the nodal points of eigenfunctions for the problem (1)-(5). We point out
that our results are extension and/or generalization to those in [7-11,
17-18, 28, 30, 31]. For example, if the retardation function $\Delta \equiv
0 $ in (1) and $p(x)\equiv 1,\delta =1,\gamma =1$ we have the formula of the
first regularized trace for the classical Sturm-Liouville operator which is
called Gelfand-Levitan formula (see [29]).

%%%%%%%%%%%%%%%%%%%%%%%%%%%%%%%%%%%%%%%%%%%%%%%%%%%%%%%%%%%%%%%%%%%

\section{The spectrum}

%%%%%%%%%%%%%%%%%%%%%%%%%%%%%%%%%%%%%%%%%%%%%%%%%%%%%%%%%%%%%%%%%%%

Let $\omega _{1}(x,\lambda )$ be a solution of Eq. (1) on $\left[ 0,\frac{%
\pi }{2}\right] ,$ satisfying the initial conditions%
\begin{equation}
\omega _{1}\left( 0,\lambda \right) =a_{2}\text{ \textit{and} }\omega
_{1}^{\prime }\left( 0,\lambda \right) =-a_{1}\text{.}  \tag{6}
\end{equation}%
The conditions (6) define a unique solution of Eq. (1) on $\left[ 0,\frac{%
\pi }{2}\right] $ (see [7,17]).

After defining the above solution, then we will define the solution $\omega
_{2}\left( x,\lambda \right) $ of Eq. (1) on $\left[ \frac{\pi }{2},\pi %
\right] $ by means of the solution $\omega _{1}\left( x,\lambda \right) $
using the initial conditions%
\begin{equation}
\omega _{2}\left( \frac{\pi }{2},\lambda \right) =\gamma _{1}\delta
_{1}^{-1}\omega _{1}\left( \frac{\pi }{2},\lambda \right) \text{ \textit{and}%
}\quad \omega _{2}^{\prime }(\frac{\pi }{2},\>\lambda )=\gamma _{2}\delta
_{2}^{-1}\omega _{1}^{\prime }(\frac{\pi }{2},\lambda )\text{.}  \tag{7}
\end{equation}%
The conditions (7) define a unique solution of Eq. (1) on $\left[ \frac{\pi 
}{2},\pi \right] $ (see [17]).

Consequently, the function $\omega \left( x,\lambda \right) $ is defined on $%
\left[ 0,\frac{\pi }{2}\right) \cup \left( \frac{\pi }{2},\pi \right] $ by
the equality%
\begin{equation*}
\omega (x,\lambda )=\left\{ 
\begin{array}{cc}
\omega _{1}(x,\lambda ), & x\in \lbrack 0,\frac{\pi }{2}), \\ 
\omega _{2}(x,\lambda ), & x\in \left( \frac{\pi }{2},\pi \right]%
\end{array}%
\right.
\end{equation*}%
is a solution of (1) on $\left[ 0,\frac{\pi }{2}\right) \cup \left( \frac{%
\pi }{2},\pi \right] $, which satisfies one of the boundary conditions and
the interface conditions (4)-(5) Then the following integral equations
hold:\ 
\begin{align}
\omega _{1}(x,\lambda )& =a_{2}\cos \frac{\lambda }{p_{1}}x-\frac{a_{1}p_{1}%
}{\lambda }\sin \frac{\lambda }{p_{1}}x  \notag \\
& -\frac{1}{\lambda p_{1}}\int\limits_{0}^{{x}}q\left( \tau \right) \sin 
\frac{\lambda }{p_{1}}\left( x-\tau \right) \omega _{1}\left( \tau -\Delta
\left( \tau \right) ,\lambda \right) d\tau ,  \tag{8}
\end{align}%
\begin{align}
\omega _{2}(x,\lambda )& =\frac{\gamma _{1}}{\delta _{1}}\omega _{1}\left( 
\frac{\pi }{2},\lambda \right) \cos \frac{\lambda }{p_{2}}\left( x-\frac{\pi 
}{2}\right) +\frac{\gamma _{2}p_{2}\omega _{1}^{\prime }\left( \frac{\pi }{2}%
,\lambda \right) }{\lambda \delta _{2}}\sin \frac{\lambda }{p_{2}}\left( x-%
\frac{\pi }{2}\right)  \notag \\
& -\frac{1}{\lambda p_{2}}\int\limits_{\frac{\pi }{2}}^{{x}}q\left( \tau
\right) \sin \frac{\lambda }{p_{2}}\left( x-\tau \right) \omega _{2}\left(
\tau -\Delta \left( \tau \right) ,\lambda \right) d\tau .  \tag{9}
\end{align}%
Solving the equations (8)-(9) by the method of successive approximation, we
obtain the following asymptotic equalities for $\left\vert \lambda
\right\vert \rightarrow \infty :$%
\begin{equation*}
\omega _{1}(x,\lambda )=a_{2}\cos \frac{\lambda }{p_{1}}x-\frac{a_{1}p_{1}}{%
\lambda }\sin \frac{\lambda }{p_{1}}x-\frac{a_{2}}{2\lambda p_{1}}%
\int\limits_{0}^{x}q(\tau )\sin \frac{\lambda }{p_{1}}\left( x-\Delta \left(
\tau \right) \right) d\tau
\end{equation*}%
\begin{equation}
-\frac{a_{2}}{2\lambda p_{1}}\int\limits_{0}^{x}q(\tau )\sin \frac{\lambda }{%
p_{1}}\left( x-\left( 2\tau -\Delta \left( \tau \right) \right) \right)
d\tau +O\left( \frac{1}{\lambda ^{2}}\right) .  \tag{10}
\end{equation}%
Differentiating (10) with respect to $x$, we get%
\begin{equation*}
\omega _{1}^{\prime }(x,\lambda )=-\frac{a_{2}\lambda }{p_{1}}\sin \frac{%
\lambda }{p_{1}}x-a_{1}\cos \frac{\lambda }{p_{1}}x-\frac{a_{2}}{2p_{1}^{2}}%
\int\limits_{0}^{x}q(\tau )\cos \frac{\lambda }{p_{1}}\left( x-\Delta \left(
\tau \right) \right) d\tau
\end{equation*}%
\begin{equation}
-\frac{a_{2}}{2p_{1}^{2}}\int\limits_{0}^{x}q(\tau )\cos \frac{\lambda }{%
p_{1}}\left( x-\left( 2\tau -\Delta \left( \tau \right) \right) \right)
d\tau +O\left( \frac{1}{\lambda }\right) .  \tag{11}
\end{equation}%
Using the fact that,%
\begin{equation*}
O\left( \frac{1}{\lambda }\right) =\left\{ 
\begin{array}{c}
\int\limits_{0}^{x}q(\tau )\sin \frac{\lambda }{p_{1}}\left( 2\tau -\Delta
\left( \tau \right) \right) d\tau ,\text{ \ }x\in \left[ 0,\frac{\pi }{2}%
\right] ; \\ 
\int\limits_{0}^{x}q(\tau )\cos \frac{\lambda }{p_{1}}\left( 2\tau -\Delta
\left( \tau \right) \right) d\tau ,\text{ \ }x\in \left[ 0,\frac{\pi }{2}%
\right] ; \\ 
\int\limits_{\frac{\pi }{2}}^{x}q(\tau )\sin \frac{\lambda }{p_{2}}\left(
2\tau -\Delta \left( \tau \right) \right) d\tau ,\text{ \ }x\in \left[ \frac{%
\pi }{2},\pi \right] ; \\ 
\int\limits_{\frac{\pi }{2}}^{x}q(\tau )\cos \frac{\lambda }{p_{2}}\left(
2\tau -\Delta \left( \tau \right) \right) d\tau ,\text{ \ }x\in \left[ \frac{%
\pi }{2},\pi \right]%
\end{array}%
\right.
\end{equation*}%
(see [9]), (10) and (11) we have%
\begin{equation*}
\omega _{2}(x,\lambda )=\frac{\gamma _{1}a_{2}}{\delta _{1}}\cos \frac{%
\lambda }{p_{2}}\left( \frac{\pi \left( p_{2}-p_{1}\right) }{2p_{1}}%
+x\right) -\frac{\gamma _{1}a_{1}p_{1}}{\lambda \delta _{1}}\sin \frac{%
\lambda }{p_{2}}\left( \frac{\pi \left( p_{2}-p_{1}\right) }{2p_{1}}+x\right)
\end{equation*}%
\begin{equation*}
-\frac{\gamma _{1}a_{2}\left( p_{1}+p_{2}\right) }{2\lambda \delta
_{1}p_{1}p_{2}}\left\{ \left[ B\left( \frac{\pi }{2},\lambda \right)
+D(x,\lambda )\right] \sin \frac{\lambda }{p_{2}}\left( \frac{\pi \left(
p_{2}-p_{1}\right) }{2p_{1}}+x\right) \right.
\end{equation*}%
\begin{equation}
\left. +\left[ A\left( \frac{\pi }{2},\lambda \right) +C(x,\lambda )\right]
\cos \frac{\lambda }{p_{2}}\left( \frac{\pi \left( p_{2}-p_{1}\right) }{%
2p_{1}}+x\right) \right\} +O\left( \frac{1}{\lambda ^{2}}\right) .  \tag{12}
\end{equation}%
Here,%
\begin{eqnarray*}
A\left( x,\lambda \right) &=&\int\limits_{0}^{x}q(\tau )\sin \frac{\lambda
\Delta \left( \tau \right) }{p_{1}}d\tau \text{ \ }x\in \left[ 0,\frac{\pi }{%
2}\right] ; \\
B\left( x,\lambda \right) &=&\int\limits_{0}^{x}q(\tau )\cos \frac{\lambda
\Delta \left( \tau \right) }{p_{1}}d\tau \text{ \ }x\in \left[ 0,\frac{\pi }{%
2}\right] ; \\
C\left( x,\lambda \right) &=&\int\limits_{\frac{\pi }{2}}^{x}q(\tau )\sin 
\frac{\lambda \Delta \left( \tau \right) }{p_{2}}d\tau ,\text{ \ }x\in \left[
\frac{\pi }{2},\pi \right] ; \\
D\left( x,\lambda \right) &=&\int\limits_{\frac{\pi }{2}}^{x}q(\tau )\cos 
\frac{\lambda \Delta \left( \tau \right) }{p_{2}}d\tau ,\text{ \ }x\in \left[
\frac{\pi }{2},\pi \right] .
\end{eqnarray*}

Differentiating (12) with respect to $x$, we get%
\begin{equation}
\omega _{2}^{\prime }(x,\lambda )=-\frac{\gamma _{1}a_{2}\lambda }{\delta
_{1}p_{2}}\sin \frac{\lambda }{p_{2}}\left( \frac{\pi \left(
p_{2}-p_{1}\right) }{2p_{1}}+x\right) -\frac{\gamma _{1}a_{1}p_{1}}{\delta
_{1}p_{2}}\cos \frac{\lambda }{p_{2}}\left( \frac{\pi \left(
p_{2}-p_{1}\right) }{2p_{1}}+x\right)  \notag
\end{equation}%
\begin{equation*}
-\frac{\gamma _{1}a_{2}\left( p_{1}+p_{2}\right) }{2\delta _{1}p_{1}p_{2}^{2}%
}\left\{ \left[ B\left( \frac{\pi }{2},\lambda \right) +D(x,\lambda )\right]
\cos \frac{\lambda }{p_{2}}\left( \frac{\pi \left( p_{2}-p_{1}\right) }{%
2p_{1}}+x\right) \right.
\end{equation*}%
\begin{equation}
\left. -\left[ A\left( \frac{\pi }{2},\lambda \right) +C(x,\lambda )\right]
\sin \frac{\lambda }{p_{2}}\left( \frac{\pi \left( p_{2}-p_{1}\right) }{%
2p_{1}}+x\right) \right\} +O\left( \frac{1}{\lambda }\right) .  \tag{13}
\end{equation}%
The solution $\omega (x,\lambda )$ defined above is a nontrivial solution of
(1) satisfying conditions (2) and (4)-(5). Putting$\>\omega (x,\lambda )\>$%
into (3), we get the characteristic equation 
\begin{equation}
\Theta (\lambda )\equiv \omega ^{\prime }(\pi ,\lambda )+d\omega (\pi
,\lambda ).  \tag{14}
\end{equation}%
The set of eigenvalues of boundary value problem (1)-(5) coincides with the
set of the squares of roots of (14), and eigenvalues are simple (see [17]).
From (12), (13) and (14), we obtain%
\begin{equation*}
\Theta (\lambda )\equiv -\frac{\gamma _{1}a_{2}\lambda }{\delta _{1}p_{2}}%
\sin \frac{\lambda \pi }{p_{2}}\left( \frac{p_{2}-p_{1}}{2p_{1}}+1\right) -%
\frac{\gamma _{1}a_{1}p_{1}}{\delta _{1}p_{2}}\cos \frac{\lambda \pi }{p_{2}}%
\left( \frac{p_{2}-p_{1}}{2p_{1}}+1\right)
\end{equation*}%
\begin{equation*}
-\frac{\gamma _{1}a_{2}\left( p_{1}+p_{2}\right) }{2\delta _{1}p_{1}p_{2}^{2}%
}\left\{ \left[ B\left( \frac{\pi }{2},\lambda \right) +D(\pi ,\lambda )%
\right] \cos \frac{\lambda \pi }{p_{2}}\left( \frac{p_{2}-p_{1}}{2p_{1}}%
+1\right) \right.
\end{equation*}%
\begin{equation*}
\left. -\left[ A\left( \frac{\pi }{2},\lambda \right) +C(\pi ,\lambda )%
\right] \sin \frac{\lambda \pi }{p_{2}}\left( \frac{p_{2}-p_{1}}{2p_{1}}%
+1\right) \right\}
\end{equation*}%
\begin{equation*}
+\frac{d\gamma _{1}a_{2}}{\delta _{1}}\cos \frac{\lambda \pi }{p_{2}}\left( 
\frac{p_{2}-p_{1}}{2p_{1}}+1\right) +O\left( \frac{1}{\lambda }\right) .
\end{equation*}

Define%
\begin{equation*}
\Theta _{0}(\lambda )\equiv -\frac{\gamma _{1}a_{2}\lambda }{\delta _{1}p_{2}%
}\sin \frac{\lambda \pi }{p_{2}}\left( \frac{p_{2}-p_{1}}{2p_{1}}+1\right) .
\end{equation*}%
Denote by $\lambda _{n}^{0}=\frac{2p_{1}p_{2}n}{p_{1}+p_{2}},n\in 
%TCIMACRO{\U{2124} }%
%BeginExpansion
\mathbb{Z}
%EndExpansion
,$ zeros of the function $\Theta _{0}\left( \lambda \right) $. It is simple
algebraically except for $\lambda _{\pm 0}^{0}$ and we have

\begin{equation}
\lambda _{n}\sim \lambda _{n}^{0}+O\left( \frac{1}{n}\right) .  \tag{15}
\end{equation}

Denote by $C_{n}$ the circle of radius, $0<\varepsilon <\frac{1}{2},$
centered at the origin $\lambda _{n}^{0}$ and by $\Gamma _{N_{0}}$ the
counterclockwise square contours with four vertices%
\begin{eqnarray*}
K &=&N_{0}+\varepsilon +N_{0}i,\text{ \ \ }L=-N_{0}-\varepsilon +N_{0}i, \\
M &=&-N_{0}-\varepsilon -N_{0}i,,\text{ \ \ }N=N_{0}+\varepsilon -N_{0}i,
\end{eqnarray*}%
where $i=\sqrt{-1}$ and $N_{0}$ is a natural number. Obviously, if $\lambda
\in C_{n}$ or $\lambda \in \Gamma _{N_{0}},$ then $\left\vert \Theta
_{0}\left( \lambda \right) \right\vert \geq M\left\vert \lambda \right\vert
e^{\left\vert \func{Im}\lambda \right\vert \pi }$ $\left( M>0\right) $ by
using a similar method in [11, 32]. Thus, on $\lambda \in C_{n}$ or $\lambda
\in \Gamma _{N_{0}},$ we have%
\begin{equation*}
\frac{\Theta \left( \lambda \right) }{\Theta _{0}(\lambda )}=1+\frac{1}{%
\lambda }\left\{ \frac{a_{1}p_{1}}{a_{2}}\cot \frac{\lambda \pi }{p_{2}}%
\left( \frac{p_{2}-p_{1}}{2p_{1}}+1\right) \right.
\end{equation*}%
\begin{equation*}
+\frac{p_{1}+p_{2}}{2p_{1}p_{2}}\left[ \left( B\left( \frac{\pi }{2},\lambda
\right) +D(\pi ,\lambda )\right) \cot \frac{\lambda \pi }{p_{2}}\left( \frac{%
p_{2}-p_{1}}{2p_{1}}+1\right) -\left( A\left( \frac{\pi }{2},\lambda \right)
+C(\pi ,\lambda )\right) \right]
\end{equation*}%
\begin{equation*}
\left. -dp_{2}\cot \frac{\lambda \pi }{p_{2}}\left( \frac{p_{2}-p_{1}}{2p_{1}%
}+1\right) \right\} +O\left( \frac{1}{\lambda ^{2}}\right) .
\end{equation*}%
Expanding $\ln \frac{\Theta \left( \lambda \right) }{\Theta _{0}\left(
\lambda \right) }$ by the Maclaurin formula, we find that%
\begin{equation*}
\ln \frac{\Theta \left( \lambda \right) }{\Theta _{0}(\lambda )}=\frac{1}{%
\lambda }\left\{ \left[ \frac{a_{1}p_{1}}{a_{2}}+\frac{p_{1}+p_{2}}{%
2p_{1}p_{2}}\left( B\left( \frac{\pi }{2},\lambda \right) +D(\pi ,\lambda
)\right) -dp_{2}\right] \right.
\end{equation*}%
\begin{equation*}
\left. \times \cot \frac{\lambda \pi }{p_{2}}\left( \frac{p_{2}-p_{1}}{2p_{1}%
}+1\right) -\frac{p_{1}+p_{2}}{2p_{1}p_{2}}\left( A\left( \frac{\pi }{2}%
,\lambda \right) +C(\pi ,\lambda )\right) \right\}
\end{equation*}%
\begin{equation*}
-\frac{1}{2\lambda ^{2}}\left\{ \left[ \frac{a_{1}p_{1}}{a_{2}}+\frac{%
p_{1}+p_{2}}{2p_{1}p_{2}}\left( B\left( \frac{\pi }{2},\lambda \right)
+D(\pi ,\lambda )\right) -dp_{2}\right] ^{2}\right.
\end{equation*}%
\begin{equation*}
\times \cot ^{2}\frac{\lambda \pi }{p_{2}}\left( \frac{p_{2}-p_{1}}{2p_{1}}%
+1\right) +\frac{\left( p_{1}+p_{2}\right) ^{2}}{4p_{1}^{2}p_{2}^{2}}\left(
A\left( \frac{\pi }{2},\lambda \right) +C(\pi ,\lambda )\right) ^{2}
\end{equation*}%
\begin{equation*}
+\frac{p_{1}+p_{2}}{p_{1}p_{2}}\left[ \frac{a_{1}p_{1}}{a_{2}}+\frac{%
p_{1}+p_{2}}{2p_{1}p_{2}}\left( B\left( \frac{\pi }{2},\lambda \right)
+D(\pi ,\lambda )\right) -dp_{2}\right]
\end{equation*}%
\begin{equation*}
\left. \times \left( A\left( \frac{\pi }{2},\lambda \right) +C(\pi ,\lambda
)\right) \cot \frac{\lambda \pi }{p_{2}}\left( \frac{p_{2}-p_{1}}{2p_{1}}%
+1\right) \right\} +O\left( \frac{1}{\lambda ^{3}}\right) .
\end{equation*}%
Using the well-known Rouche Theorem, we get that $\Theta \left( \lambda
\right) $ has the same number of zeros inside $\Gamma _{N_{0}}$ as $\Theta
_{0}\left( \lambda \right) $ (see [11]). Using the residue theorem, we have%
\begin{equation*}
\lambda _{n}-\lambda _{n}^{0}=-\frac{1}{2\pi i}\doint\limits_{C_{n}}\ln 
\frac{\Theta \left( \lambda \right) }{\Theta _{0}\left( \lambda \right) }%
d\lambda
\end{equation*}%
\begin{equation*}
=-\frac{1}{2\pi i}\doint\limits_{C_{n}}\left[ \frac{a_{1}p_{1}}{a_{2}}+\frac{%
p_{1}+p_{2}}{2p_{1}p_{2}}\left( B\left( \frac{\pi }{2},\lambda \right)
+D(\pi ,\lambda )\right) -dp_{2}\right] \frac{\cot \frac{\lambda \pi }{p_{2}}%
\left( \frac{p_{2}-p_{1}}{2p_{1}}+1\right) }{\lambda }d\lambda
\end{equation*}%
\begin{equation*}
+\frac{1}{2\pi i}\doint\limits_{C_{n}}\frac{p_{1}+p_{2}}{2p_{1}p_{2}}\left(
A\left( \frac{\pi }{2},\lambda \right) +C(\pi ,\lambda )\right) \frac{%
d\lambda }{\lambda }
\end{equation*}%
\begin{equation*}
+\frac{1}{2\pi i}\doint\limits_{C_{n}}\left[ \frac{a_{1}p_{1}}{a_{2}}+\frac{%
p_{1}+p_{2}}{2p_{1}p_{2}}\left( B\left( \frac{\pi }{2},\lambda \right)
+D(\pi ,\lambda )\right) -dp_{2}\right] ^{2}\frac{\cot ^{2}\frac{\lambda \pi 
}{p_{2}}\left( \frac{p_{2}-p_{1}}{2p_{1}}+1\right) }{2\lambda ^{2}}d\lambda
\end{equation*}%
\begin{equation*}
+\frac{1}{2\pi i}\doint\limits_{C_{n}}\frac{\left( p_{1}+p_{2}\right) ^{2}}{%
4p_{1}^{2}p_{2}^{2}}\left( A\left( \frac{\pi }{2},\lambda \right) +C(\pi
,\lambda )\right) ^{2}\frac{d\lambda }{2\lambda ^{2}}
\end{equation*}%
\begin{equation*}
+\frac{1}{2\pi i}\doint\limits_{C_{n}}\frac{p_{1}+p_{2}}{p_{1}p_{2}}\left[ 
\frac{a_{1}p_{1}}{a_{2}}+\frac{p_{1}+p_{2}}{2p_{1}p_{2}}\left( B\left( \frac{%
\pi }{2},\lambda \right) +D(\pi ,\lambda )\right) -dp_{2}\right]
\end{equation*}%
\begin{equation*}
\times \left( A\left( \frac{\pi }{2},\lambda \right) +C(\pi ,\lambda
)\right) \frac{\cot \frac{\lambda \pi }{p_{2}}\left( \frac{p_{2}-p_{1}}{%
2p_{1}}+1\right) }{2\lambda ^{2}}d\lambda +O\left( \frac{1}{n^{3}}\right) .
\end{equation*}%
Thus, using (15) and residue calculation we have proven the following
theorem.

\begin{theorem}
The spectrum of the problem (1)-(5) has the%
\begin{equation*}
\lambda _{n}=\lambda _{n}^{0}-\frac{1}{\lambda _{n}^{0}\pi }\left[ \frac{%
a_{1}p_{1}}{a_{2}}+\frac{p_{1}+p_{2}}{2p_{1}p_{2}}\left( B\left( \frac{\pi }{%
2},\lambda _{n}^{0}\right) +D(\pi ,\lambda _{n}^{0})\right) -dp_{2}\right]
\end{equation*}%
\begin{equation*}
+\frac{1}{\left( \lambda _{n}^{0}\right) ^{2}\pi }\left\{ \frac{p_{1}+p_{2}}{%
p_{1}p_{2}}\left( A\left( \frac{\pi }{2},\lambda _{n}^{0}\right) +C(\pi
,\lambda _{n}^{0})\right) \right.
\end{equation*}%
\begin{equation*}
\left. \times \left[ \frac{a_{1}p_{1}}{a_{2}}+\frac{p_{1}+p_{2}}{2p_{1}p_{2}}%
\left( B\left( \frac{\pi }{2},\lambda _{n}^{0}\right) +D(\pi ,\lambda
_{n}^{0})\right) -dp_{2}\right] \right\}
\end{equation*}%
\begin{equation*}
-\frac{1}{\left( \lambda _{n}^{0}\right) ^{3}}\left[ \frac{a_{1}p_{1}}{a_{2}}%
+\frac{p_{1}+p_{2}}{2p_{1}p_{2}}\left( B\left( \frac{\pi }{2},\lambda
_{n}^{0}\right) +D(\pi ,\lambda _{n}^{0})\right) -dp_{2}\right] ^{2}+O\left( 
\frac{1}{n^{3}}\right)
\end{equation*}%
asymptotic distribution for sufficiently large $\left\vert n\right\vert .$
\end{theorem}

\section{The regularized trace formula}

In this section, we will get regularized trace formula for the problem
(1)-(5).

The asymptotic formula (15) for the eigenvalues implies that for all
sufficiently large $N_{0},$ the numbers $\lambda _{n}$ with $\left\vert
n\right\vert \leq N_{0}$ are inside $\Gamma _{N_{0}},$ and the numbers $%
\lambda _{n}$ with $\left\vert n\right\vert >N_{0}$ are outside $\Gamma
_{N_{0}}.$ It follows that%
\begin{equation*}
\lambda _{-0}^{2}+\lambda _{0}^{2}+\sum_{0\neq n=-N_{0}}^{N_{0}}\left(
\lambda _{n}^{2}-\left( \lambda _{n}^{0}\right) ^{2}\right) =-\frac{1}{2\pi i%
}\doint\limits_{\Gamma _{n}}2\lambda \ln \frac{\Delta \left( \lambda \right) 
}{\Delta _{0}\left( \lambda \right) }d\lambda
\end{equation*}%
\begin{equation*}
=-\frac{1}{2\pi i}\doint\limits_{\Gamma _{n}}2\left[ \frac{a_{1}p_{1}}{a_{2}}%
+\frac{p_{1}+p_{2}}{2p_{1}p_{2}}\left( B\left( \frac{\pi }{2},\lambda
\right) +D(\pi ,\lambda )\right) -dp_{2}\right]
\end{equation*}%
\begin{equation*}
\times \cot \frac{\lambda \pi }{p_{2}}\left( \frac{p_{2}-p_{1}}{2p_{1}}%
+1\right) d\lambda +\frac{1}{2\pi i}\doint\limits_{\Gamma _{n}}2\left[ \frac{%
p_{1}+p_{2}}{2p_{1}p_{2}}\left( A\left( \frac{\pi }{2},\lambda \right)
+C(\pi ,\lambda )\right) \right] d\lambda
\end{equation*}%
\begin{equation*}
+\frac{1}{2\pi i}\doint\limits_{\Gamma _{n}}\left[ \frac{a_{1}p_{1}}{a_{2}}+%
\frac{p_{1}+p_{2}}{2p_{1}p_{2}}\left( B\left( \frac{\pi }{2},\lambda \right)
+D(\pi ,\lambda )\right) -dp_{2}\right] ^{2}\frac{\cot ^{2}\frac{\lambda \pi 
}{p_{2}}\left( \frac{p_{2}-p_{1}}{2p_{1}}+1\right) }{\lambda }d\lambda
\end{equation*}%
\begin{equation*}
+\frac{1}{2\pi i}\doint\limits_{\Gamma _{n}}\frac{\left( p_{1}+p_{2}\right)
^{2}}{4p_{1}^{2}p_{2}^{2}}\left( A\left( \frac{\pi }{2},\lambda \right)
+C(\pi ,\lambda )\right) ^{2}\frac{d\lambda }{\lambda }
\end{equation*}%
\begin{equation*}
+\frac{1}{2\pi i}\doint\limits_{\Gamma _{n}}\frac{p_{1}+p_{2}}{p_{1}p_{2}}%
\left[ \frac{a_{1}p_{1}}{a_{2}}+\frac{p_{1}+p_{2}}{2p_{1}p_{2}}\left(
B\left( \frac{\pi }{2},\lambda \right) +D(\pi ,\lambda )\right) -dp_{2}%
\right]
\end{equation*}%
\begin{equation*}
\times \left( A\left( \frac{\pi }{2},\lambda \right) +C(\pi ,\lambda
)\right) \frac{\cot \frac{\lambda \pi }{p_{2}}\left( \frac{p_{2}-p_{1}}{%
2p_{1}}+1\right) }{\lambda }d\lambda +O\left( \frac{1}{N_{0}}\right) ,
\end{equation*}%
by calculations, which implies that%
\begin{equation*}
\lambda _{-0}^{2}+\lambda _{0}^{2}+\sum_{0\neq n=-N_{0}}^{N_{0}}\lambda
_{n}^{2}-\left( \lambda _{n}^{0}\right) ^{2}
\end{equation*}%
\begin{equation*}
=-\frac{2}{\pi }\sum_{0\neq n=-N_{0}}^{N_{0}}\left[ \frac{a_{1}p_{1}}{a_{2}}+%
\frac{p_{1}+p_{2}}{2p_{1}p_{2}}\left( B\left( \frac{\pi }{2},\lambda
_{n}^{0}\right) +D(\pi ,\lambda _{n}^{0})\right) -dp_{2}\right]
\end{equation*}%
\begin{equation*}
-\frac{2}{\pi }\left[ \frac{a_{1}p_{1}}{a_{2}}+\frac{p_{1}+p_{2}}{2p_{1}p_{2}%
}\left( B\left( \frac{\pi }{2},0\right) +D(\pi ,0)\right) -dp_{2}\right]
\end{equation*}%
\begin{equation*}
+\sum_{0\neq n=-N_{0}}^{N_{0}}\frac{p_{1}+p_{2}}{p_{1}p_{2}}\left[ \frac{%
a_{1}p_{1}}{a_{2}}+\frac{p_{1}+p_{2}}{2p_{1}p_{2}}\left( B\left( \frac{\pi }{%
2},\lambda _{n}^{0}\right) +D(\pi ,\lambda _{n}^{0})\right) -dp_{2}\right]
\end{equation*}%
\begin{equation*}
\times \left( A\left( \frac{\pi }{2},\lambda _{n}^{0}\right) +C(\pi ,\lambda
_{n}^{0})\right) \frac{1}{\lambda _{n}^{0}\pi }+R
\end{equation*}%
\begin{equation*}
-\left[ \frac{a_{1}p_{1}}{a_{2}}+\frac{p_{1}+p_{2}}{2p_{1}p_{2}}\left(
B\left( \frac{\pi }{2},0\right) +D(\pi ,0)\right) -dp_{2}\right] ^{2}
\end{equation*}%
\begin{equation}
+\frac{\left( p_{1}+p_{2}\right) ^{2}}{4p_{1}^{2}p_{2}^{2}}\left( A\left( 
\frac{\pi }{2},0\right) +C(\pi ,0)\right) ^{2}+O\left( \frac{1}{N_{0}}%
\right) ,  \tag{16}
\end{equation}%
where,%
\begin{equation*}
R=Res_{\lambda =0}\left\{ \frac{p_{1}+p_{2}}{p_{1}p_{2}}\left[ \frac{%
a_{1}p_{1}}{a_{2}}+\frac{p_{1}+p_{2}}{2p_{1}p_{2}}\left( B\left( \frac{\pi }{%
2},\lambda \right) +D(\pi ,\lambda )\right) -dp_{2}\right] \right.
\end{equation*}%
\begin{equation*}
\left. \times \left( A\left( \frac{\pi }{2},\lambda \right) +C(\pi ,\lambda
)\right) \frac{\cot \frac{\lambda \pi }{p_{2}}\left( \frac{p_{2}-p_{1}}{%
2p_{1}}+1\right) }{\lambda }d\lambda \right\} .
\end{equation*}%
Passing to the limit as $N_{0}\rightarrow \infty $ in (16), we have%
\begin{equation*}
\lambda _{-0}^{2}+\lambda _{0}^{2}+\sum_{0\neq n=-\infty }^{\infty }\left\{
\lambda _{n}^{2}-\left( \lambda _{n}^{0}\right) ^{2}\right.
\end{equation*}%
\begin{equation*}
+\frac{2}{\pi }\left[ \frac{a_{1}p_{1}}{a_{2}}+\frac{p_{1}+p_{2}}{2p_{1}p_{2}%
}\left( B\left( \frac{\pi }{2},\lambda _{n}^{0}\right) +D(\pi ,\lambda
_{n}^{0})\right) -dp_{2}\right]
\end{equation*}%
\begin{equation*}
-\frac{p_{1}+p_{2}}{p_{1}p_{2}}\left[ \frac{a_{1}p_{1}}{a_{2}}+\frac{%
p_{1}+p_{2}}{2p_{1}p_{2}}\left( B\left( \frac{\pi }{2},\lambda
_{n}^{0}\right) +D(\pi ,\lambda _{n}^{0})\right) -dp_{2}\right]
\end{equation*}%
\begin{equation*}
\left. \times \left( A\left( \frac{\pi }{2},\lambda _{n}^{0}\right) +C(\pi
,\lambda _{n}^{0})\right) \frac{1}{\lambda _{n}^{0}\pi }\right\}
\end{equation*}%
\begin{equation*}
=-\frac{2}{\pi }\left[ \frac{a_{1}p_{1}}{a_{2}}+\frac{p_{1}+p_{2}}{%
2p_{1}p_{2}}\left( B\left( \frac{\pi }{2},0\right) +D(\pi ,0)\right) -dp_{2}%
\right] +R
\end{equation*}%
\begin{equation*}
-\left[ \frac{a_{1}p_{1}}{a_{2}}+\frac{p_{1}+p_{2}}{2p_{1}p_{2}}\left(
B\left( \frac{\pi }{2},0\right) +D(\pi ,0)\right) -dp_{2}\right] ^{2}
\end{equation*}%
\begin{equation}
+\frac{\left( p_{1}+p_{2}\right) ^{2}}{4p_{1}^{2}p_{2}^{2}}\left( A\left( 
\frac{\pi }{2},0\right) +C(\pi ,0)\right) ^{2}.  \tag{17}
\end{equation}%
The series on the left side of (17) is called the regularized trace of the
problem (1)-(5).

\section{The oscillation}

In this chapter, we will find nodal points of eigenfunctions of the problem
(1)-(5).

Let us rewrite the equation (10) and replace $\lambda $ by $\lambda _{n}$%
\begin{equation*}
\omega _{1}(x,\lambda _{n})=a_{2}\cos \frac{\lambda _{n}x}{p_{1}}-\frac{%
a_{1}p_{1}}{\lambda _{n}}\sin \frac{\lambda _{n}x}{p_{1}}-\frac{a_{2}\sin 
\frac{\lambda _{n}x}{p_{1}}}{2\lambda _{n}p_{1}}\int\limits_{0}^{x}q(\tau
)\cos \left( \frac{\lambda _{n}\Delta \left( \tau \right) }{p_{1}}\right)
d\tau
\end{equation*}%
\begin{equation*}
+\frac{a_{2}\cos \frac{\lambda _{n}x}{p_{1}}}{2\lambda _{n}p_{1}}%
\int\limits_{0}^{x}q(\tau )\sin \left( \frac{\lambda _{n}\Delta \left( \tau
\right) }{p_{1}}\right) d\tau +O\left( \frac{1}{\lambda _{n}^{2}}\right) .
\end{equation*}%
Let us assume that $x_{n}^{j}$ are the nodal points of the eigenfunction $%
\omega _{1}\left( x,\lambda _{n}\right) .$ Taking $\sin \left( \frac{\lambda
_{n}x}{p_{1}}\right) \neq 0$ into account for sufficiently large $n,$ we get 
\begin{equation*}
\cot \left( \frac{\lambda _{n}x}{p_{1}}\right) \left[ 1+\frac{A\left(
x,\lambda _{n}\right) }{2\lambda _{n}p_{1}}\right] =\frac{a_{1}p_{1}}{%
a_{2}\lambda _{n}}+\frac{B\left( x,\lambda _{n}\right) }{2\lambda _{n}p_{1}}%
+O\left( \frac{1}{\lambda _{n}^{2}}\right) .
\end{equation*}%
It follows easily that 
\begin{equation}
\tan \left( \frac{\lambda _{n}x}{p_{1}}+\frac{\pi }{2}\right) =-\frac{%
a_{1}p_{1}}{a_{2}\lambda _{n}}-\frac{B\left( x,\lambda _{n}\right) }{%
2\lambda _{n}p_{1}}+O\left( \frac{1}{\lambda _{n}^{2}}\right) .  \tag{18}
\end{equation}%
Thus, solving the equation (18), one obtains%
\begin{equation}
x_{n}^{j}=\frac{\left( j-\frac{1}{2}\right) \pi p_{1}}{\lambda _{n}}-\frac{%
a_{1}p_{1}^{2}}{a_{2}\lambda _{n}^{2}}-\frac{B\left( x_{n}^{j},\lambda
_{n}\right) }{2\lambda _{n}^{2}}+O\left( \frac{1}{\lambda _{n}^{3}}\right) .
\tag{19}
\end{equation}%
Note that%
\begin{equation}
\lambda _{n}^{-1}=\frac{1}{\lambda _{n}^{0}}+\frac{\left[ \frac{a_{1}p_{1}}{%
a_{2}}+\frac{p_{1}+p_{2}}{2p_{1}p_{2}}\left( B\left( \frac{\pi }{2},\lambda
_{n}^{0}\right) +D(\pi ,\lambda _{n}^{0})\right) -dp_{2}\right] }{\left(
\lambda _{n}^{0}\right) ^{3}\pi }+O\left( \frac{1}{n^{4}}\right)  \tag{20}
\end{equation}%
and%
\begin{equation}
\lambda _{n}^{-2}=\frac{1}{\left( \lambda _{n}^{0}\right) ^{2}}+O\left( 
\frac{1}{n^{4}}\right) .  \tag{21}
\end{equation}%
Substituting (20) and (21) into (19) we have 
\begin{equation*}
x_{n}^{j}=\frac{\left( j-\frac{1}{2}\right) \pi p_{1}}{\lambda _{n}^{0}}-%
\frac{\left( j-\frac{1}{2}\right) p_{1}\left[ \frac{a_{1}p_{1}}{a_{2}}+\frac{%
p_{1}+p_{2}}{2p_{1}p_{2}}\left( B\left( \frac{\pi }{2},\lambda
_{n}^{0}\right) +D(\pi ,\lambda _{n}^{0})\right) -dp_{2}\right] }{\left(
\lambda _{n}^{0}\right) ^{3}}
\end{equation*}%
\begin{equation}
-\frac{a_{1}p_{1}^{2}}{a_{2}\left( \lambda _{n}^{0}\right) ^{2}}-\frac{%
B\left( \frac{\left( j-\frac{1}{2}\right) \pi p_{1}}{\lambda _{n}^{0}}%
,\lambda _{n}^{0}\right) }{2\left( \lambda _{n}^{0}\right) ^{2}}+O\left( 
\frac{1}{n^{3}}\right) ,\text{ \ \ \ }j=\overline{1,\left[ \frac{n}{2}\right]
}.  \tag{22}
\end{equation}%
Similarly, from (12), we get 
\begin{equation*}
\frac{\delta _{1}\omega _{2}(x,\lambda _{n})}{\gamma _{1}}=a_{2}\cos \frac{%
\lambda _{n}}{p_{2}}\left( \frac{\pi \left( p_{2}-p_{1}\right) }{2p_{1}}%
+x\right) -\frac{a_{1}p_{1}}{\lambda _{n}}\sin \frac{\lambda _{n}}{p_{2}}%
\left( \frac{\pi \left( p_{2}-p_{1}\right) }{2p_{1}}+x\right)
\end{equation*}%
\begin{equation*}
-\frac{a_{2}\left( p_{1}+p_{2}\right) }{2\lambda _{n}p_{1}p_{2}}\left\{ %
\left[ B\left( \frac{\pi }{2},\lambda _{n}\right) +D(x,\lambda _{n})\right]
\sin \frac{\lambda _{n}}{p_{2}}\left( \frac{\pi \left( p_{2}-p_{1}\right) }{%
2p_{1}}+x\right) \right.
\end{equation*}%
\begin{equation*}
\left. +\left[ A\left( \frac{\pi }{2},\lambda _{n}\right) +C(x,\lambda _{n})%
\right] \cos \frac{\lambda _{n}}{p_{2}}\left( \frac{\pi \left(
p_{2}-p_{1}\right) }{2p_{1}}+x\right) \right\} +O\left( \frac{1}{\lambda
_{n}^{2}}\right) .
\end{equation*}%
For nodal points of $\omega _{2}(x,\lambda _{n})$, again, taking $\sin \frac{%
\lambda _{n}}{p_{2}}\left( \frac{\pi \left( p_{2}-p_{1}\right) }{2p_{1}}%
+x\right) \neq 0$ into account for sufficiently large $n,$ we get%
\begin{eqnarray*}
&&\cot \frac{\lambda _{n}}{p_{2}}\left( \frac{\pi \left( p_{2}-p_{1}\right) 
}{2p_{1}}+x\right) \left[ 1-\frac{\left( A\left( \frac{\pi }{2},\lambda
_{n}\right) +C\left( x,\lambda _{n}\right) \right) \left( p_{1}+p_{2}\right) 
}{2\lambda _{n}p_{1}p_{2}}\right] \\
&=&\frac{a_{1}p_{1}}{a_{2}\lambda _{n}}+\frac{\left( p_{1}+p_{2}\right)
\left( B\left( \frac{\pi }{2},\lambda _{n}\right) +D\left( x,\lambda
_{n}\right) \right) }{2\lambda _{n}p_{1}p_{2}}+O\left( \frac{1}{\lambda
_{n}^{2}}\right) .
\end{eqnarray*}%
and thus%
\begin{eqnarray}
&&\tan \left( \frac{\lambda _{n}}{p_{2}}\left( \frac{\pi \left(
p_{2}-p_{1}\right) }{2p_{1}}+x\right) +\frac{\pi }{2}\right)  \notag \\
&=&-\frac{a_{1}p_{1}}{a_{2}\lambda _{n}}-\frac{\left( p_{1}+p_{2}\right)
\left( B\left( \frac{\pi }{2},\lambda _{n}\right) +D\left( x,\lambda
_{n}\right) \right) }{2\lambda _{n}p_{1}p_{2}}+O\left( \frac{1}{\lambda
_{n}^{2}}\right) .  \TCItag{23}
\end{eqnarray}%
Thus, solving the equation (23), one obtains%
\begin{equation*}
x_{n}^{j}=-\frac{\pi \left( p_{2}-p_{1}\right) }{2p_{1}}+\frac{\left( j-%
\frac{1}{2}\right) \pi p_{2}}{\lambda _{n}}-\frac{\alpha _{1}p_{1}p_{2}}{%
a_{2}\lambda _{n}^{2}}
\end{equation*}%
\begin{equation}
-\frac{\left( p_{1}+p_{2}\right) \left( B\left( \frac{\pi }{2},\lambda
_{n}\right) +D\left( x_{n}^{j},\lambda _{n}\right) \right) }{2\lambda
_{n}^{2}p_{1}}+O\left( \frac{1}{\lambda _{n}^{3}}\right) .  \tag{24}
\end{equation}%
Substituting (20) and (21) into (24) we have 
\begin{equation*}
x_{n}^{j}=-\frac{\pi \left( p_{2}-p_{1}\right) }{2p_{1}}+\frac{\left( j-%
\frac{1}{2}\right) \pi p_{2}}{\lambda _{n}^{0}}
\end{equation*}%
\begin{equation*}
-\frac{\left( j-\frac{1}{2}\right) p_{2}\left[ \frac{a_{1}p_{1}}{a_{2}}+%
\frac{p_{1}+p_{2}}{2p_{1}p_{2}}\left( B\left( \frac{\pi }{2},\lambda
_{n}^{0}\right) +D(\pi ,\lambda _{n}^{0})\right) -dp_{2}\right] }{\left(
\lambda _{n}^{0}\right) ^{3}}-\frac{\alpha _{1}p_{1}p_{2}}{a_{2}\left(
\lambda _{n}^{0}\right) ^{2}}
\end{equation*}%
\begin{equation}
-\frac{\left( p_{1}+p_{2}\right) \left[ B\left( \frac{\pi }{2},\lambda
_{n}^{0}\right) +D\left( \frac{\left( j-\frac{1}{2}\right) \pi p_{2}}{%
\lambda _{n}^{0}},\lambda _{n}^{0})\right) \right] }{2\left( \lambda
_{n}^{0}\right) ^{2}p_{1}}+O\left( \frac{1}{n^{3}}\right) ,\text{ \ \ \ }j=%
\overline{\left[ \frac{n}{2}\right] +1,n}.  \tag{25}
\end{equation}%
Thus we have proven the following theorem:

\begin{theorem}
For sufficiently large $n$, we have the \ formulas (22) and (25) of the
nodal points for the problem (1)-(5).
\end{theorem}

%%%%%%%%%%%%%%%%%%%%%%%%%%%%%%%%%%%%%%%%%%%%%%%%%%%%%%%%%%%%%%%%


\begin{thebibliography}{99}
\bibitem{} L. Crocco and S. Chang, Theory of combustion instability in
liquid propellant rocket motors, Butterworths, London, 1956.

\bibitem{} J. Kolesovas, D. Svitra, \emph{Mathematical modelling of the
combustion process in the chamber of a liquid propellant rocket engine},
Lithuanian Mathematical Journal, 15 (4) (1975) 632-642.

\bibitem{} K.F. Teodorcik, Self-oscillatory systems, 3rd ed., Gostehizdat,
Moscow, 1952. (Russian).

\bibitem{} A.A. Harkevic, Auto-oscillations, Gostehizdat, Moscow, 1954.
(Russian).

\bibitem{} A. V. Likov, Y. A. Mikhalilov, The Theory of Heat and Mass
Transfer, Qosenergaizdat, 1963 (In Russian).

\bibitem{} I. Titeux, Y. Yakubov, \emph{Completeness of root functions for
thermal conduction in a strip with piecewise continuous coefficients}, Math.
Models Methods Appl. Sci. 7 (7) (1997) 1035-1050.

\bibitem{} S.B. Norkin, Differential equations of the second order with
retarded argument, Translations of Mathematical Monographs, AMS, Providence,
RI, 1972.

\bibitem{} M. Pikula, \emph{Regularized traces of differential operator of
Sturm-Liouville type with retarded argument}, Differentsialnye Uravneniya,
26 (1) (1990) 103-109 (Russian); English translation: Differential
Equations, 26 (1) (1990) 92-96.

\bibitem{} A. Bayramov, S. \d{C}al\i \d{s}kan and S. Uslu, \emph{Computation
of eigenvalues and eigenfunctions of a discontinuous boundary value problem
with retarded argument}, Appl. Math. Comput., 191 (2007) 592-600.

\bibitem{} E. \c{S}en, A. Bayramov, \emph{Calculation of eigenvalues and
eigenfunctions of a discontinuous boundary value problem with retarded
argument which contains a spectral parameter in the boundary condition},
Math. Comput. Model. 54 (11-12) (2011) 3090-3097.

\bibitem{} C-F. Yang, \emph{Trace and inverse problem of a discontinuous
Sturm-Liouville operator with retarded argument}, J. Math. Anal. Appl., 395
(2012) 30-41.

\bibitem{} F.A. Cetinkaya, K.R. Mamedov, \emph{A boundary value problem with
retarded argument and discontinuous coefficient in the differential equation}%
, Azerbaijan Journal of Mathematics, 7 (2017) 135-145.

\bibitem{} E. \c{S}en, M. Acikgoz, S. Araci, \emph{Spectral problem for
Sturm-Liouville operator with retarded argument which contains a spectral
parameter in the boundary condition}, Ukrainian Mathematical Journal, 68 (8)
(2017) 1263-1277.

\bibitem{} F. Hira, \emph{A trace formula for the Sturm-Liouville type
equation with retarded argument}, Commun. Fac. Sci. Univ. Ank. Ser. A1 Math.
Stat., 66 (1) (2017) 124-132.

\bibitem{} M. Bayramoglu, A. Bayramov, E. \c{S}en, \emph{A regularized trace
formula for a discontinuous Sturm-Liouville operator with delayed argument},
Electronic Journal of Differential Equations (EJDE), 2017 (104) (2017) 1-12.

\bibitem{} E. \c{S}en, \emph{A regularized trace formula and oscillation of
eigenfunctions of a Sturm-Liouville operator with retarded argument at two
points of discontinuity},\emph{\ }Mathematical Methods in the Applied
Sciences, doi: 10.1002/mma.4512.

\bibitem{} E. \c{S}en, A. Bayramov, \emph{Spectral analysis of boundary
value problems with retarded argument}, Commun. Fac. Sci. Univ. Ank. S\'{e}%
r. A1 Math. Stat., 66 (2) (2017) 175-194.

\bibitem{} I.M. Gelfand, B.M. Levitan, \emph{On a formula for eigenvalues of
a differential operator of second order}, Doklady Akademii Nauk SSSR, 88
(1953) 593--596 (Russian).

\bibitem{} A. Makin, \emph{Regularized trace of the Sturm-Liouville operator
with irregular boundary conditions,} Electronic Journal of Differential
Equations (EJDE), 2009 (27) (2009) 1-8.

\bibitem{} A. Makin, \emph{Trace formulas for the Sturm-Liouville operator
with regular boundary conditions}, Doklady Mathematics, 76 (2) (2007)
702-707.

\bibitem{} C-F. Yang, \emph{Trace formula for the matrix Sturm-Liouville
operator}, Anal. Math. Phys. 6(2016) 31--41.

\bibitem{} F.G. Maksudov, M. Bayramoglu, E.E. Ad\i guzelov, \emph{On a
regularized traces of the Sturm-Liouville operator on a finite interval with
the unbounded operator coefficient}, Doklady Akademii Nauk SSSR, 277 (4)
(1984); English translation: Soviet Math. Dokl., 30 (1984) 169-173.

\bibitem{} E. Adiguzelov, Y. Sezer, \emph{The regularized trace of a self
adjoint differential operator of higher order with unbounded operator
coefficient}, Appl. Math. Comput., 218 (2011) 2113-2121.

\bibitem{} E. \c{S}en, A. Bayramov, K. Oru\c{c}o\u{g}lu; \emph{Regularized
trace formula for higher order differential operators with unbounded
coefficients}, Electronic Journal of Differential Equations (EJDE), 2016
(31) (2016), 1-12.

\bibitem{} E. Adiguzelov, Y. Sezer, \emph{The second regularized trace of a
self adjoint differential operator given in a finite interval with bounded
operator coefficient}, Math. Comput. Model., 53 (2011) 553-565.

\bibitem{} F. Gesztesy F., H. Holden, On Trace Formulas for Schr\"{o}%
dinger-Type Operators. In: Truhlar D.G., Simon B. (eds) Multiparticle
Quantum Scattering With Applications to Nuclear, Atomic and Molecular
Physics. The IMA Volumes in Mathematics and its Applications, vol 89.
Springer, New York, 1997.

\bibitem{} V.A. Sadovnichii, V.E. Podolskii, \emph{Traces of operators},
Uspekhi Mat. Nauk, 61 (2006) 89-156; English translation: Russian Math.
Surveys 61 (2006) 885-953.

\bibitem{} C.T. Fulton and S.A. Pruess, \emph{Eigenvalue and eigenfunction
asymptotics for regular Sturm-Liouville problems}, J. Math. Anal. Appl. 188
(1994) 297-340.

\bibitem{} L.A. Dikii, \emph{Trace formulas for Sturm-Liouville differential
equations}, Uspekhi Mat. Nauk, 12 (3) (1958) 111-143.

\bibitem{} C. T. Shieh, V.A. Yurko, \emph{Inverse nodal and inverse spectral
problems for discontinuous boundary value problem}, J. Math. Anal. Appl. 347
(2008) 266--272.

\bibitem{} Y.H. Cheng, C.K. Law, \emph{On the quasi-nodal map for the
Sturm-Liouville problem}, Proc. Soc. Edinburgh, 136A (2006) 71-86.

\bibitem{} V.A. Yurko, Inverse Spectral Problems for Differential Operators
and Their Applications, Gordon and Breach, Amsterdam, 2000.

%%%%%%%%%%%%%%%%%%%%%%%%%%%%%%%%%%%%%%%%%%%%%%%%%%%%%%%%%%%%%%%%%%%
\end{thebibliography}
\end{document}